\documentclass[a4paper, 11pt, leqno]{article}
\usepackage{amsmath, amsfonts, amssymb}
\usepackage{graphics}
\usepackage{color}
\usepackage{fancyhdr} 
\usepackage{sectsty}
\usepackage[T1]{fontenc}
\usepackage{pstricks}
\usepackage{mathrsfs}
\usepackage{pifont}
\usepackage{array}
\usepackage{titlesec}

\usepackage{ae}
\usepackage{aecompl}

\newcommand{\real}{\mathbb{R}}
\newcommand{\integer}{\mathbb{Z}}

\newcommand{\restr}[1]{|_{#1}}
\newcommand{\htop}{h_{\text{top}}}

\newcommand{\s}{\mathscr S}
\newcommand{\il}{\mathscr L}
\newcommand{\F}{\mathcal F}
\newcommand{\I}{\mathscr I}

\newcommand{\diff}[1]{\text{Diff}\,(#1)}
\newcommand{\diffrel}[2]{\text{Diff}\,(#1, \text{rel } #2)}

\newcommand{\diffrelid}[3]{\text{Diff}\,(#1, \text{rel } #2, #3)}
\newcommand{\iso}{\sim}
\newcommand{\br}[1]{\text{Br}\,(#1)}

\newcommand{\card}{\#}

\titleformat{\section}{}{\Large \bf \thesection.}{0.3em}{\bfseries\Large}
\titlespacing{\section}{0cm}{1.cm}{0.4cm}

\makeatletter
\@addtoreset{figure}{section}
\makeatother

\newcounter{theolistcounter}
\newenvironment{theolist}{%
\begin{list}{(\roman{theolistcounter})}{%
\usecounter{theolistcounter}
\setlength{\itemsep}{0cm}}}%
{\end{list}}

\newcounter{alglistcounter}
{\end{list}}

\newcounter{defn}[section]
\renewcommand{\thedefn}{\arabic{section}.\arabic{defn}}
\newenvironment{defn}[0]%
{\refstepcounter{defn}\vspace{10pt}\par\noindent 
{\bf Definition \thedefn} \ -- \ 
\begin{itshape}}%
{\end{itshape}\par\vspace{0.2cm}}

\newcounter{prop}[section]
\renewcommand{\theprop}{\arabic{section}.\arabic{prop}}
\newenvironment{prop}[0]%
{\refstepcounter{prop}\vspace{10pt}\par\noindent 
{\bf Proposition \theprop} \ -- \ 
\begin{itshape}}%
{\end{itshape}\par\vspace{0.2cm}}

\newcounter{theo}[section]
\renewcommand{\thetheo}{\arabic{section}.\arabic{theo}}
\newenvironment{theo}[0]%
{\refstepcounter{theo}\vspace{10pt}\par\noindent 
{\bf Theoreme \thetheo} \ -- \
\begin{itshape}}%
{\end{itshape}\par\vspace{0.2cm}} 

\newcounter{ex}[section]
\renewcommand{\theex}{\arabic{section}.\arabic{ex}}
\newenvironment{ex}[0]%
{\refstepcounter{ex}\vspace{10pt}\par\noindent 
{\bf Example \theex} \ -- \ }%
{\par\vspace{0.2cm}} 

\newcounter{conj}[section]
\renewcommand{\theconj}{\arabic{section}.\arabic{conj}}
\newenvironment{conj}[0]%
{\refstepcounter{conj}\vspace{10pt}\par\noindent 
{\bf Conjecture \theconj} \ -- \ 
\begin{itshape}}%
{\end{itshape}\par\vspace{0.2cm}}

\newenvironment{dem}{\noindent {\it Proof. }} %
{$\square$ \par \vspace{0.2cm}}

\begin{document}

\begin{center}
{\Large  On Computing the Entropy of Braids} \\[0.2cm]
{\it Jacques-Olivier Moussafir}\\[0.4cm]
{\small Saint-Gobain Recherche \\
39 quai Lucien Lefranc \\
93300 Aubervillers, France}
\end{center}

\vspace{0.4cm}
\noindent
{\bf Abstract.} 
We consider in this paper the problem of computing the 
entropy of a braid. We recall its definition and construct, 
for each braid, a sequence of real numbers, whose limit is 
the braid's entropy.
We state one conjecture about the convergence speed, 
and two about the braids that have high entropy, 
but are written with few letters.

\footnotetext[1]{Mathematics Subject Classification: 37B40, 37E30}


\section{Introduction}
The formal definition of braid's entropy,
as introduced by P. Boyland in \cite{Boyland}, is a very
natural extension of results presented in \cite{FLP} about 
dynamical properties of pseudo-Anosov maps. 

In this paper we present a method for computing the entropy of braids
without working with train tracks, see \cite{BH}. More precisely, we
construct for every braid $\beta$, an integer sequence $c_m$
whose growth factor is the entropy of $\beta$. 
The major drawback of this method is that it is not an algorithm 
since we don't get the result in finitely many steps: 
we have to use an artificial stopping criterion. 
The major advantage over train-tracks is that it works faster. 

The paper is organized as follows. In the first section we quote
some results from \cite{FLP} that justify Boyland's definition and
will be used as technical components of proofs. Next section is the
central technical part of the paper. Informally, we connect Boyland's
definition with Dynnikov and Wiest braid complexity, see \cite{CB}. 
More precisely, we construct for each surface diffeomorphism
$\varphi$, an integer sequence $c_m$ such that 
$$
\lim_{m\rightarrow \infty} \frac{1}{m} \log c_m = \htop (\varphi),
$$
where $\htop (\varphi)$ stands for the topological entropy of
$\varphi$, see \cite{KH} for definitions and properties. The following
section gives the formal definition of the braid entropy, introduces
integral laminations and their coding as they appear in \cite{YB}, and
finally presents a formula that gives the minimum number of
intersections of an integral lamination with the real axis in terms of
Dynnikov's coordinates. At this point we describe the method for
the computation of braid entropies, say a few words about the
corresponding computer program and conjecture it might be turned 
into an algorithm. Last section presents some conjectures about the braids 
that are written with few letters but have large entropy.

\section{Some results on surface diffeomorphisms}
One can find in \cite{FLP} a complete proof of Thurston's results
about the classification of surface diffeomorphisms. But the authors also
prove many propositons and lemmas that help to understand Thurston's
result and that we shall need. Let us quote some of them.

Let  $M$ be a compact oriented surface possibly with boundary, then 
$\s(M)$ will denote the set of homotopy classes of closed and connected 
simple paths that are not homotopic to zero or to a component of 
the boundary of $M$.

If $(\F, \mu)$ is a measured foliation of $M$, see \cite{FLP} for 
a definition, and $\gamma$ a closed and connected simple path then
$$
\int_\gamma |\mu|
$$
will denote the total variation of the differential form $\mu$ along
$\gamma$. If $\alpha \in \s(M)$, let
$$
\I(\F,\mu,\alpha) = \inf_{\gamma \in \alpha} \int_{\gamma} |\mu|.
$$

\begin{prop}
  Let $\varphi$ be a pseudo-Anosov map with stable and unstable
  foliations $(\F^s, \mu^s)$ and $(\F^u, \mu^u)$. Let $\alpha \in
  \s(M)$, then 
  $$
  \I(\F^s,\mu^s,\alpha)>0 \quad \text{and} 
  \quad \I(\F^u,\mu^u,\alpha) >0.
  $$
\end{prop}

For any two $\alpha, \beta \in \s(M)$, $c(\alpha, \beta)$ will denote
the minimum number intersections of $a$ and $b$, with $a \in \alpha$
and $b \in \beta$. 

Let $M$ be a compact oriented surface, $A$ and $B$ two subsets of $M$.
We'll note $\diff{M}$ the group of diffeomorphisms of $M$,
$\diffrel{M}{A}$ for the subgroup of $\diff{M}$ whose elements
leave $A$ invariant, and finally $\diffrelid{M}{A}{B}$ for the subgroup of
$\diffrel{M}{A}$ whose elements are the identity on $B$.

If $G$ is a subgroup of $\diff{M}$ (such as $\diffrel{M}{A}$),
$\varphi$, $\psi \in G$, we'll write $\varphi \iso_G \psi$ or
simply $\varphi \iso \psi$ when $\varphi$ and $\psi$ are isotopic in $G$.

\begin{prop}\label{prop htop_pa}
  Let $\varphi$ be a pseudo-Anosov map with stable and unstable
  foliations $(\F^s, \mu^s)$, $(\F^u, \mu^u)$ and let $\lambda$ be a
  real number, $\lambda>1$, be such that $\varphi(\F^s) = 1/\lambda\, 
  \F^s$ and $\varphi(\F^u) = \lambda \, \F^u$.
  Then
  \begin{equation}
  \lim_{n\rightarrow \infty} \frac{c(\varphi^n \alpha,
    \beta)}{\lambda^n} = \I(\F^s,\mu^s,\alpha) \,
    \I(\F^u,\mu^u,\alpha), 
  \end{equation}
  \begin{equation}
    \htop(\varphi) = \log \lambda,
  \end{equation}
  \begin{equation}
  \htop(\varphi) = \inf \{ 
  \htop(\psi), \ \psi \in \diff{M}, \ \psi \iso \varphi \}.
  \end{equation}

\end{prop}
The following result is a formulation of Thurston's result about the
classification of surface diffeomorphisms.

\begin{theo}\label{theo thurston}
  Let  $M$ be a compact oriented surface possibly with boundary, and
  $\varphi$ a diffeomorphism of $M$. 
  There exists finitely many simple curves of $M$: $C_1, \ldots, C_l$ 
  and a diffeomorphism $\psi$ isotopic to $\varphi$ such that 
  cutting $M$ along $C_1, \ldots, C_l$ produces $M_1, \ldots, M_k$
  compact, possibly non connected, oriented surfaces 
  satisfying
  \begin{theolist}
  \item $M = M_1 \cup \ldots \cup M_k$;
  \item if $j_1\neq j_2$, $M_{j_1} \cap M_{j_2}$ is the union 
    of several $C_i$;
  \item $\varphi$ is isotopic to $\psi$;
  \item for $i=1,\ldots, l$, $\psi(M_i) = M_i$;
  \item for $i=1,\ldots, l$, $\psi\restr{M_i}$ is periodic or pseudo-Anosov.
  \end{theolist} 
\end{theo}

\begin{ex}
  Let $M=\real^2 / \integer^2$ be the two-dimensional torus, and
  $\varphi$ be the Dehn twist given by
  $$
  \begin{pmatrix}
    1 & 1 \\
    0 & 1
  \end{pmatrix}.
  $$
  In this case the simple closed curve $C_1 = [(0,0), (1,0)]$
  decomposes $M$ into $C_1$ and a closed cylinder, $M_1$.
  The map $\varphi\restr{M_1}$ is isotopic to the identity.
\end{ex}

\section{Counting intersections}
Proposition \ref{prop htop_pa} already suggests a way of computing the
topological entropy of a pseudo-Anosov maps $\varphi$ of a compact
oriented surface $M$: pick $\alpha$ and $\beta$ in $\s (M)$, and
compute 
$$
\lim_{m\rightarrow \infty} \frac{1}{m} \log c(\varphi^m \alpha,
\beta).$$
To do so, one would have to choose some $\alpha$ and $\beta$, and
devise a method for the computation of $c(\varphi^m \alpha, \beta)$.

\begin{figure}[hbt!]
  \begin{center}
    \input{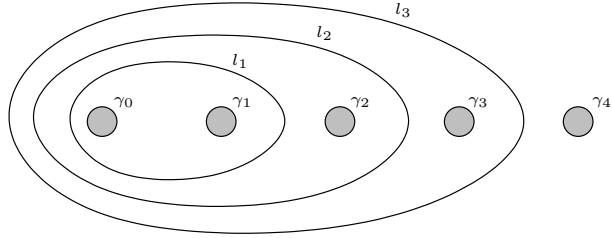}
    \caption{\label{fig L} A sphere with $n=6$ open discs removed 
      is represented on the plane using a stereographic projection.
      The disc $\gamma_5$ containing the projection's pole is omitted.
      In this case the set of curves $L$ is made of $3$ 
      disjoint simple closed curves: $l_1$, $l_2$, $l_3$.}
  \end{center}
\end{figure}

\begin{figure}[hbt!]
  \begin{center}
    \input{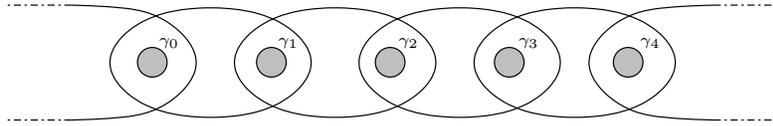}
    \caption{\label{fig R} The set $R$ corresponding to a sphere with six
      holes. Here it contains 6 simple closed curves.}
  \end{center}
\end{figure}

To go on with this idea, we first restrict ourselves to surfaces $M$
that are homeomorphic to a sphere with finitely many open discs removed.
The boundaries of removed open discs will usually be denoted 
$\gamma_0$, \ldots, $\gamma_n$, see figure 
\ref{fig L} or \ref{fig R}.
We'll usually represent such surfaces on the plane using the
stereographic projection with respect to a pole belonging to one of
the removed discs. The boundary of the disc containing the pole 
will usually not be drawn.

We'll consider in $M$ two sets of curves, up to homotopy. 
The first set $L$ is represented on figure \ref{fig L}, 
and the second set, $R$, on figure \ref{fig R}. 

If $A$ and $B$ are two finite subsets of $\s(M)$, $c(A,B)$ will denote
$$
c(A,B) = \sum_{a_i \in A, \, b_j \in B} c(a_i, b_j).
$$

\begin{prop}\label{prop main_result}
  Let $M$ be a sphere with $n \ge 3$ discs removed,
  let $\varphi$ be a diffeomorphism of $M$ and let $L$ and $R$
  be the two subsets of $\s(M)$ represented on figure \ref{fig L} and
  \ref{fig R}.
  Then,
  $$
  \lim_{m\rightarrow \infty} \frac{1}{m} \log c(\varphi^m L, R) = 
  \htop (\varphi).
  $$
\end{prop}

\begin{dem}
  Let $\gamma_0$,\ldots, $\gamma_n$ denote the boundaries of the
  removed discs, and
  let $\psi$, $M_1$, \ldots, $M_k$ be the diffeomorphism and 
  the surfaces given by \ref{theo thurston}. 
  Then, let
  $$
  I_j = \left\{ i \in \{0,\ldots, n \}, \ \gamma_i \subset M_j \right\}.
  $$
 
  The decomposition $M=M_1 \cup \ldots M_k$ induces a decomposition of
  $l\in \s(M)$ into finitely many pieces -- we consider sufficiently
  regular elements in $l$ and assume $M_j$ have smooth boundaries with
  finitely many singularities.

  For each component of $l$ that belongs to $M_j$, we form a new simple 
  closed curve $l'$ with that component and connect its ends following the
  boundary of $M_j$. There are two choices for this. If there are more
  than three $\gamma_i$ inside $M_j$ we choose $l'$ such that $l' \in
  \s(M_j)$. We shall denote $l_{I_j}$ the union of all these curves, see
  figure \ref{fig lij}.
  \begin{figure}[hbt!]
    \begin{center}
      \input{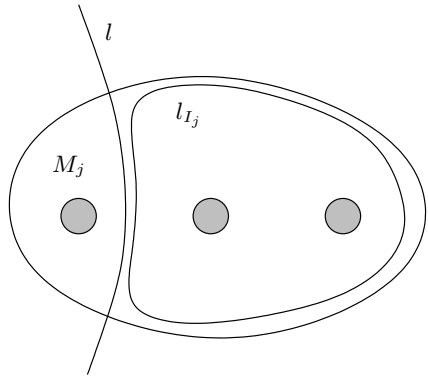}
      \caption{\label{fig lij} Construction of $l_{I_j}$ from $l$. If there
	are more than three holes in $M_j$, we construct $l_{I_j}$ in such a
	way that it belongs to $\s(M_j)$.}
    \end{center}
  \end{figure}

  We'll also form two sets of curves $R_{I_j}$ and $R_{I_j}'$ from $R$
  and $M_j$. First, for $i,j \in {0, \ldots n}$, let $r_{i,j}$ be
  the simple closed curve formed with the boundaries of $\gamma_i$ and
  $\gamma_j$. Figure \ref{fig R} represents for instance, from left to
  right, $r_{5,0}$, $r_{0,1}$, $r_{1,2}$, $r_{2,3}$, $r_{3,4}$ and
  $r_{4,5}$. Let 
  $$
  I_j = \{ i_1, i_2, \ldots, i_p \}, \quad \text{with} \quad i_1 < i_2 < \ldots
  < i_p,
  $$
  then, 
  $$
  R_{I_j} = r_{i_1,i_2} \cup r_{i_2,i_3} \ldots \cup r_{i_{p-1},i_p}
  \cup r_{i_p,i_1} 
  $$
  and
  $$
  R_{I_j}' = \bigcup_{
    \begin{array}{l}
      \gamma_i \text{ or } \gamma_{i+1} \text{ is inside }\\
      M_j \text{ but not both.}
    \end{array}
  } 
  r_{i, i+1}.
  $$
  This construction is represented on figure \ref{fig rij}.
  \begin{figure}[hbt!]
    \begin{center}
      \input{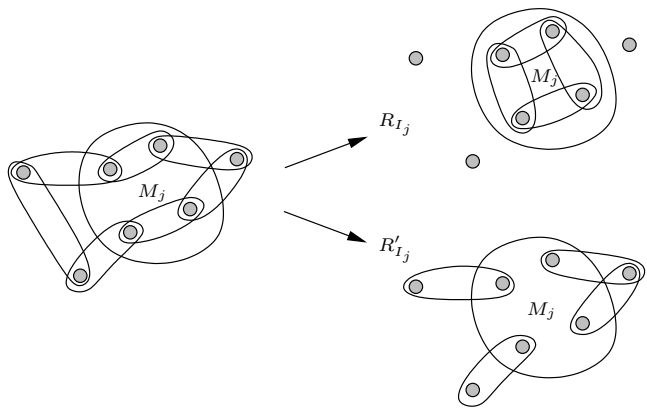}
      \caption{\label{fig rij} Construction of $R_{I_j}$ and $R_{I_j}'$
	from $R$ and $M_j$. Each $\gamma_i$ is represented.}
    \end{center}
  \end{figure}
  For any $l \in \s(M)$, $c(l, R)$ may be written as follows
  $$
  c(l, R) = \sum_{
    \begin{array}{l}
      l_{I_j} \in \s(M_j) \\ 
      \card I_j \ge 2
  \end{array}}
  c(l_{I_j}, R_{I_j})
  \ \ + \ \
  \sum_j c(l, R_{I_j}'),
  $$
  where $\card I_j$ is the cardinal of $I_j$.
  \begin{figure}[hbt!]
    \begin{center}
      \input{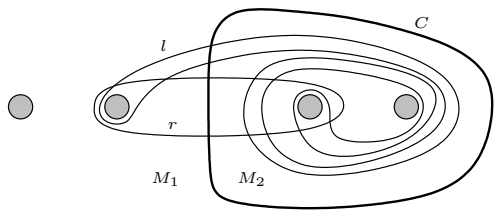}
      \caption{\label{fig dehn}A Dehn twist $\delta$ in $M_2$ along $C$ 
	produces a twist of $l$, and $c(\delta^m l, r)$ grows linearly with $m$.}
    \end{center}
  \end{figure}
  For any $l \in \s(M)$, the sequence $c(\psi^m l, R_{I_j}')$ cannot
  grow faster than linearly with $m$ because only Dehn twists along
  $M_j$'s boundaries can make it grow, as illustrated on figure \ref{fig
    dehn}, hence
  $$
  \lim_{m\rightarrow \infty} \frac{1}{m} \log c(\psi^m l, R) = 
  \max_{l_{I_j} \in \s(M_j),\ \card I_j \ge 2} 
  \ 
  \lim_{m\rightarrow \infty} \frac{1}{m} \log c(\varphi^m l_{I_j}, R_{I_j}). 
  $$
  Remark that if $\card {I_j}\le 2$, then 
  $$
  \lim \frac{1}{m} \log c(\varphi^m l_{I_j}, R_{I_j}) = 0.
  $$
  Let $h_j$ be the topological entropy
  of $\psi\restr{M_j}$, using proposition \ref{prop htop_pa} we get
  $$
  \lim_{m\rightarrow \infty} \frac{1}{m} \log c(\psi^m l, R) = 
  \max_{l_{I_j} \in \s(M_j), \ \card I_j \ge 2} \ h_j.
  $$
  If $\card I_j \ge 3$, there exists $l\in L$ such that
  $l_{I_j} \in \s(M_j)$, hence
  $$
  \begin{array}{ll}
  \lim_{m\rightarrow \infty} \frac{1}{m} \log c(\psi^m L, R) 
  & = \max_{\card I_j \ge 3} \, h_j \\[0.2cm]
  & = \max_{j\in \{1,\ldots, k\} } \, h_j
  \end{array}
  $$
  This last maximum is precisely the topological entropy of $\psi$,
  and consequently, the topological entropy of $\varphi$.
\end{dem}

\section{The entropy of braids}
\subsection{Definitions}
The braid group can be defined in many ways, see \cite{Birman} for
instance. Here, we shall use a definition that shows the connection
with surface diffeomorphism.

\begin{defn}  
  Let $n$ be an integer, $n\ge 2$, and $M$ the surface obtained from the unit
  disc $D^2$ with $n$ disjoint open discs removed; 
  we'll write $\Gamma_n$ the union of their boundaries.
  The braid group with $n$ strands is defined as 
  $$
  \br{n} = \diffrelid{D^2}{\Gamma_n}{\partial D^2} / \iso
  $$
\end{defn}

\begin{defn}
  Let $\beta \in \br{n}$. The entropy of $\beta$
  is defined as
  $$
  h(\beta) = \inf_{\varphi \in \beta} \htop (\varphi).
  $$
\end{defn}
Proposition \ref{prop htop_pa} should make this definition rather clear.
\subsection{Integral laminations}
We already came across integral laminations in proposition \ref{prop
  main_result}. Here is a more formal definition.

\begin{defn}
  Let $M$  be a compact and oriented surface. An integral lamination of
  $M$ is a set $L$ of disjoint non homotopic simple closed curves
  of $M$. Integral laminations are considered up to homotopy. 
  The set of integral laminations of $M$ will be denoted 
  $\il (M)$.
\end{defn}

\begin{figure}[hbt!]
  \begin{center}
    \input{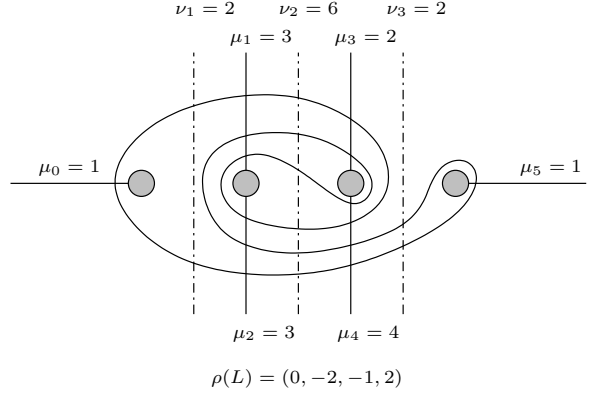}
    \caption{\label{fig lamcod} Dynnikov's coding of an integral
    lamination $L$ of $M_n$, with $n=2$. 
    There are $8$ intersections between $L$ and the
    real axis.}
  \end{center}
\end{figure}

I.~Dynnikov introduced in \cite{YB} a coding for integral laminations
when $M_n$ is a sphere with $n+3$ open discs removed.
We represented such a surface on figure \ref{fig lamcod}. 
Following our convention one component of $M$'s boudary is not represented.
We consider $2n+2$ half lines (continuous lines) 
and $n+1$ vertical lines (dotted lines)
as represented on figure \ref{fig lamcod}. We consider $\nu_1$,
\ldots, $\nu_n$ the minimum number of intersection of $L$ with each
dotted line, and $\mu_0$, $\mu_1$, $\mu_2$,\ldots, $\mu_{2n-1}$,
$\mu_{2n}$, the minimum number of intersection of $L$ with the half
lines. Now for $i=1$, \ldots $n$, let
$$
\begin{array}{l}
  a_i = \frac{\mu_{2i} - \mu_{2i - 1}}{2} \\
  b_i = \frac{\nu_{i} - \nu_{i+1}}{2}.
\end{array}
$$
and 
$$
\rho(L) = (a_1, b_1, \ldots, a_n, b_n). 
$$
The following result appears in \cite{YB}. 
\begin{prop} 
  Let $n$ be an integer, $n\ge 2$. The map $\rho$ defines a bijection 
  between $\il(M_n)$ an $\integer^{2n}$.
\end{prop}

\begin{ex}
For the integral lamination $L$ represented on figure \ref{fig L},
$\rho(L) = (0,1,0,1,0,1)$. We'll denote $L_0^n$ the integral
lamination whose coding is $(0,1,0,1,\ldots, 0,1)$.  
\end{ex}

\begin{figure}[hbt!]
  \begin{center}
    \input{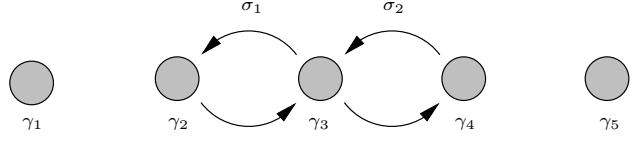}
    \caption{\label{fig brgen} Here $n=3$. The border of one disc is not
    represented. The braid group with $3$ strands acts on $M_n$. We illustrated
    the actions of the generators $\sigma_1$ and $\sigma_2$.}
  \end{center}
\end{figure}

We need to adapt slightly the definition we gave earlier 
for the braid group.
\begin{prop}
  Let $n$ be an integer, $n\ge 2$, $M_n$ the sphere with $n+3$ open discs
  removed. Let $\gamma_0$, \ldots, $\gamma_{n+2}$ denote their
  boundaries, then,
  $$
  \br{n} = 
  \diffrelid{M}{\gamma_{1} \cup \ldots \cup \gamma_{n}}
	    {\gamma_{0} \cup \gamma_{n+1} \cup \gamma_{n+2}}
	    / 
	    \iso.
	    $$
\end{prop}
Using this definition, we choose the $n-1$ generators of $\br{n}$ 
to be the diffeomorphisms whose action is depicted on \ref{fig brgen}.
One checks easily that these generators satisfy the usual relations,
see \cite{FLP} for instance.
The action of $\br{n}$ on $\il (M_n)$ is coded using I.~Dynnikov's
formulae. For any real number $a$, we'll write
$$
a^+ = \max(a,0) \quad a^-= \min(a,0).
$$

\begin{prop}
  Let $n$ be an integer $n\ge 2$, and $L \in \il (M_n)$
  with $\rho(M)=(a_1,b_1,\ldots, a_n, b_n)$.
  For each integer $i$, $1 \le i \le n-1$, let
  $$
  \begin{array}{l}
    \rho(\sigma_i L) =  (a'_1,b'_1,\ldots, a'_n, b'_n) \\[0.1cm]
    \rho(\sigma_i^{-1} L) =  (a''_1,b''_1,\ldots, a''_n, b''_n) \\[0.1cm]
    c = a_i - b_i^- -a_{i+1} + b_{i+1}^+ \\[0.1cm]
    d = a_i + b_i^- -a_{i+1} - b_{i+1}^+.
  \end{array}
  $$
  Then, 
  $$
  \begin{array}{l}
    a'_j = a_j,\ b'_j = b_j  \quad \text{for } j\neq i, i+1 \\[0.2cm]
    \left\{
    \begin{array}{l}
      a'_i = a_i + b_i^+ + (b_{i+1}^+ - c)^+\\[0.2cm] 
      b'_i = b_{i+1} - c^+\\[0.2cm] 
      a'_{i+1} = a_{i+1} + b_{i+1}^- + (b_i^- + c)^-\\[0.2cm]
      b'_{i+1} = b_i + c^+
    \end{array}
    \right.
  \end{array}
  $$
  and
  $$
  \begin{array}{l}
    a''_j = a_j,\ b''_j = b_j  \quad \text{for } j\neq i, i+1 \\[0.2cm]
    \left\{
    \begin{array}{l}
      a''_i = a_i - b_i^+ - (b_{i+1}^+ + d)^+\\[0.2cm] 
      b''_i = b_{i+1} + d^-\\[0.2cm]
      a''_{i+1} = a_{i+1} - b_{i+1}^- - (b_i^- - d)^-\\[0.2cm]
      b''_{i+1} = b_i - d^-.
    \end{array}
    \right.
  \end{array}
  $$
\end{prop}
We may use these formulas to count the minimum number of 
intersection of an integral lamination with the real axis.
\begin{prop}\label{prop inter}
Let $n$ be an integer, $n \ge 2$, $L \in \il(M_n)$ and $c(L)$
denote the minimum number of intersections between $L$ and the real
axis. If $\rho(L)=(a_1, b_1, \ldots, a_n, b_n)$, then
$$
c(L) = \sum_{i=1}^n |b_i| + \sum_{i=1}^{n-1} |a_{i+1} - a_i| 
+ |a_1| + |a_n| + \nu_1/2 + \nu_n/2.
$$
\end{prop}
\begin{dem}
Four types of intersections may occur as represented on figure 
\ref{fig lamcod}. Counting them results in the formula
for $c(L)$.
\end{dem}
If $R$ is the set of curves we introduced earlier, see \ref{fig R}, 
and $L$ the integral lamination represented on figure \ref{fig L}, then
$$
c(L,R) = 2 c(L).
$$
This leads to the description of a method for the estimation of braids
entropy. Let $n$ be an integer, $n\ge 1$, $\beta \in \br{n}$,
recall that $L_0^n$ is the lamination whose coordinate is
$(0,1,0,1,\ldots,0,1)$, and choose $\varepsilon > 0$.

\begin{dingautolist}{172}
\item Write $\beta$ as a word using standard generators
$\sigma_1$,\ldots, $\sigma_{n-1}$. 	     
$$
\beta = \prod \ \sigma_{i_k}^{a_{i_k}}, \ a_{i_k} \in \integer.
$$

\item For $m=1,2,\ldots$, compute $\rho (\beta^m L_0^n)$ and 
$c_n = 1/m \, \log c(\beta^m L_0^n)$
using Dynnikov's formulae, proposition \ref{prop inter} and forgetting
about $\nu_1/2$ and $\nu_2/2$ since they don't change.
\item Stop when $| c_{m+1} - c_m | < \varepsilon$.
\end{dingautolist}

As we said this method is not an algorithm, nevertheless the
corresponding computer program behaves well. We recover 
the well-known topological entropy of 
$\beta =\sigma_1 \sigma_2^{-1} \in \br{3}$
$$
h(\beta) = \log \frac{3 + \sqrt{5}}{2} \simeq 0.962
$$
Moreover the sequences $c_m$ seem to decrease like 
$\log m / m$. We have the following conjecture.

\begin{conj}
Let $n$ be an integer, $n\ge 2$. There exists a positive constant 
$C_n \in \real$ such that for any braid $\beta \in \br{n}$ 
and its corresponding sequence $(c_m)_{m\ge 0}$,
$$
|c_m - h(\beta)| \le C_n \frac{\log m}{m}.
$$
\end{conj}
This last conjecture, if it was correct, would turn our method into an
algorithm for the computation of $h(\beta)$ with precision $\varepsilon$.

Finally, notice that we consider integer sequences $c(\beta^m L_0^n)$
that grow rapidly, like geometric ones. 
We had to use a special library for handling large
integers. We used {\tt NTL}, see \cite{NTL}. 
Our program is freely available at 
\verb+http://www.ceremade.dauphine.fr/~msfr+.

\section{Braids with maximum entropy}
We were primarily interested in braids that have large entropy, but
are written with few letters.

\begin{defn}
Let $n$ be an integer $n\ge 2$, and $\br{n}$ the braid group with $n$
strands and standard generators $\sigma_1$, \ldots, $\sigma_{n-1}$. 
If $\beta \in \br{n}$, the length of $\beta$, $l(\beta)$ 
is the minimum number of $\sigma_i$ needed to write $\beta$.
\end{defn}
\begin{defn}
We say that a braid $\beta$ has maximal entropy if
$$
h(\beta) = \max \{ h(\beta'), \ l(\beta') \le l(\beta) \}.
$$
\end{defn}
Notice here that we do not refer to the number of strands.
Using the program we got the two following conjectures.
\begin{conj}
  Braids of maximal entropy belong to $\br{3}$ or $\br{4}$.
\end{conj}

\begin{defn}
  Let $n$ be an integer $n\ge 2$, and $\beta \in \br{n}$. If
  $$
  \beta = \sigma_{i_1}^{\epsilon_1}\ldots \sigma_{i_k}^{\epsilon_k}
  $$
  with
  $$
  i_{k+1} = i_k \pm 1 ,\ |\epsilon_k| = 1,\ \epsilon_{k+1} = -\epsilon_{k},
  $$
  then $\beta$ is said to be alternated (with respect with the
  standard generators).
\end{defn}

\begin{conj}
  Braids with maximal entropy are alternated.
\end{conj}

\vspace{0.2cm}
\noindent
{\bf Acknowledgments.} It is a pleasure to thank Jean-Luc Thiffeault
who introduced me to Boyland's ideas, Sergey Tarasov and Sergei Lando
who suggested using Dynnikov's results and Victor Turchin 
for useful conversations.
\nocite{*}
\bibliographystyle{plain}
\bibliography{EB}

\begin{thebibliography}{10}

\bibitem{ArAv}
V.I. Arnol'd and A.~Avez.
\newblock {\em Problèmes ergodiques de la mécanique classique}.
\newblock Gauthier-Villars, 1967.

\bibitem{BH}
M.~Bestvina and M.~Handel.
\newblock Train-tracks for surface homeomorphisms.
\newblock {\em Topology}, 34(1):109--140, 1995.

\bibitem{Birman}
J.~Birman.
\newblock {\em Braids, links, and mapping class groups}.
\newblock {Annals of Mathematics Studies.} {Princeton University Press and
  University of Tokyo Press.}, 1975.

\bibitem{Boyland}
P.~Boyland.
\newblock Topological methods in surface dynamics.
\newblock {\em Topology Appl.}, 58(3):223--298, 1994.

\bibitem{Boyland_Aref_Stremler}
P.~Boyland, H.~Aref, and M.~Stremler.
\newblock Topological fluid mechanics of stirring.
\newblock {\em J. Fluid Mech.}, 403:277--304, 2000.

\bibitem{YB}
I.A. Dynnikov.
\newblock {On a Yang-Baxter mapping and the Dehornoy ordering.}
\newblock {\em Russ. Math. Surv.}, 57(3):592--594, 2002.

\bibitem{CB}
I.A. Dynnikov and B.~Wiest.
\newblock On the complexity of braids.
\newblock Paper available on line at {\tt http://hal.ccsd.cnrs.fr/}, 2004.

\bibitem{FLP}
A.~Fathi, F.~Laudenbach, and V.~Poenaru.
\newblock {\em Travaux de Thurston sur les surfaces. Séminaire Orsay. 2nd ed.},
  volume 66-67 of {\em Astérisque}.
\newblock Centre National de la Recherche Scientifique, 1991.

\bibitem{Katok}
A.~Katok.
\newblock Lyapunov exponents, entropy and periodic orbits for diffeomorphisms.
\newblock {\em Publ. Math., Inst. Hautes Étud. Sci.}, 1980.

\bibitem{KH}
A.~Katok and B.~Hasselblatt.
\newblock {\em Introduction to the modern theory of dynamical systems}.
\newblock Cambridge University Press, 1995.

\bibitem{Kolev}
B.~Kolev.
\newblock {Entropie topologique et représentation de Burau}.
\newblock {\em C. R. Acad. Sci., Paris, Sér. I}, 309(13), 1989.

\bibitem{NTL}
V.~Shoup.
\newblock Ntl: A library for doing number theory.
\newblock The library is available at {\tt http://www.shoup.net}, 2001.

\bibitem{Thurston}
W.~Thurston.
\newblock On the geometry and dynamics of diffeomorphisms of surfaces.
\newblock {\em Bull. Am. Math. Soc., New Ser.}, 19(2):417--431, 1988.

\end{thebibliography}

\end{document}